# Objective Bayes testing of Poisson versus inflated Poisson models[*]

**M. J. Bayarri[1], James O. Berger[2] and Gauri S. Datta[3]**

*University of Valencia, Duke University and SAMSI, and University of Georgia*

**Abstract:** The Poisson distribution is often used as a standard model for count data. Quite often, however, such data sets are not well fit by a Poisson model because they have more zeros than are compatible with this model. For these situations, a zero-inflated Poisson (ZIP) distribution is often proposed. This article addresses testing a Poisson versus a ZIP model, using Bayesian methodology based on suitable objective priors. Specific choices of objective priors are justified and their properties investigated. The methodology is extended to include covariates in regression models. Several applications are given.

## Contents

## 1. Introduction

The Poisson distribution is often used as a standard probability model for count data. For example, a production engineer may count the number of defects in items

[*]Supported in part by NSF under Grants DMS-01-03265, SES-02-41651 and AST-05-07481, and by the Spanish Ministry of Education and Science, under Grant MTM2007-61554.
[1]Department of Statistics and O.R., University of Valencia, Av. Dr. Moliner 50, 46100 Burjassot, Valencia, Spain, e-mail: susie.bayarri@uv.es
[2]ISDS, Box 90251, Durham, NC 27708-0251, and 19 T.W. Alexander Dr., P.O. Box 14006, Research Triangle Park, NC 27709-4006, USA, e-mail: berger@samsi.info
[3]Department of Statistics, University of Georgia, Athens, GA 30602-1952, USA, e-mail: gaurisdatta@gmail.com

randomly selected from a production process. Quite often, however, such data sets are not well fit by a Poisson model because they contain more zero counts than are compatible with the Poisson model. An example is again provided by the production process; indeed, according to Ghosh et al. [14], when some production processes are in a near perfect state, zero defects will occur with a high probability. However, random changes in the manufacturing environment can lead the process to an imperfect state, producing items with defects. The production process can move randomly back and forth between the perfect and the imperfect states. For this type of production process many items will be produced with zero defects, and this excess might be better modeled by a ZIP distribution than a Poisson distribution.

For $0 \leq p \leq 1, \lambda > 0$, the ZIP$(\lambda, p)$ distribution has the probability function

$$(1.1) \qquad f_1(x \mid \lambda, p) = p\, I(x=0) + (1-p)\, f_0(x \mid \lambda), \quad x = 0, 1, 2, \ldots,$$

where $I(\cdot)$ is the indicator function, and $f_0(x|\lambda)$ is the Poisson probability function

$$(1.2) \qquad f_0(x|\lambda) = \frac{e^{-\lambda} \lambda^x}{x!}, \ x = 0, 1, 2, \ldots.$$

The parameter $p$ is referred to as the *zero-inflation parameter*.

Many authors used the ZIP distribution with and without covariates to model count data. In a ZIP regression model, Lambert [18] used a frequentist approach and Ghosh et al. [14] used a Bayesian approach to analyze industrial data sets.

While the aforementioned authors used the ZIP model to analyze their data, a number of authors have addressed the problem of checking whether a ZIP model is needed to model the data. From the frequentist perspective, score tests have been developed for testing the hypothesis $\mathcal{H}_0 : p = 0$ vs. $\mathcal{H}_1 : p \neq 0$ in a ZIP regression model ([10], [12]). From the Bayesian perspective, Bhattacharya et al. [9] presented a Bayesian method to test $p \leq 0$ versus the alternative $p > 0$ by computing a certain posterior probability of the alternative hypothesis. As in ([10], [12]), $p$ is allowed to be negative in their model [9], as long as $p + (1-p)e^{-\lambda} \geq 0$.

In this paper, we consider Bayesian testing of $M_0$ versus $M_1$ given by

$$(1.3) \qquad M_0: \quad X_i \overset{i.i.d.}{\sim} f_0(\cdot \mid \lambda), \ i = 1, \ldots, n,$$

$$(1.4) \qquad M_1: \quad X_i \overset{i.i.d.}{\sim} f_1(\cdot \mid \lambda, p), \ i = 1, \ldots, n,$$

where $f_0, f_1$ are given in (1.1) and (1.2), respectively. Note that, as opposed to the situations in the papers mentioned above, $p < 0$ is not possible here. Indeed, we can alternatively formulate the problem as that of testing, within the ZIP model,

$$\mathcal{H}_0 : p = 0 \quad \text{versus} \quad \mathcal{H}_1 : p > 0.$$

Unlike the analysis in [9], $p = 0$ (i.e., the Poisson model) is assumed to have a priori believability (e.g., prior probability 1/2).

In Section 2 we develop the suggested objective testing of Poisson versus ZIP models when not all counts are zeros. For all zeros, the ZIP distribution is not identifiable, and a proper prior is required for all parameters; we address this in Section 5. Section 3 is devoted to some comparative examples. We consider inclusion of covariates in Section 4, where we address the testing of Poisson versus ZIP regression models and give an example involving AIDS related deaths in men. In the regression case, in order for the objective Bayesian model selection to be successful we need enough positive counts so that the design matrix based on the positive counts is full column rank. When this condition does not hold we suggest in Section 5 a partially proper prior on the regression parameters to be used for model selection. Proofs and technical details are relegated to an Appendix.

## 2. Formulation of the problem

The Bayesian methodology for choosing between two models for some data is conceptually very simple (see, e.g., [3]). One assesses the prior probabilities of each model, the prior distributions for the model parameters, and computes the posterior probabilities of each model. These posterior probabilities can be computed directly from the prior probabilities and the *Bayes Factor*, an (integrated) likelihood ratio for the models which is very popular in Bayesian testing and model selection.

Often it is not possible (for lack of time or resources) to carefully assess in a subjective manner all the needed priors. In these situations, very satisfactory answers are provided by *objective Bayesian analyses* that do not use external information other than that required to formulate the problem (see [4]). First we review below some difficulties of model selection via objective Bayesian analysis. Then we justify the objective prior we chose for our problem, derive the corresponding Bayes Factor and study properties of the prior and the Bayes factor.

### 2.1. Bayesian model selection and Bayes factors

To compare two models, $M_0$ and $M_1$, for the data $\boldsymbol{X} = (X_1, \ldots, X_n)$, the Bayesian approach is based on the *Bayes factor* $B_{10}$ of $M_1$ to $M_0$ given by

$$(2.1) \qquad B_{10} = \frac{m_1(\boldsymbol{x})}{m_0(\boldsymbol{x})} = \frac{\int f_1(\boldsymbol{x} \mid \boldsymbol{\theta}_1)\pi_1(\boldsymbol{\theta}_1)d\boldsymbol{\theta}_1}{\int f_0(\boldsymbol{x} \mid \boldsymbol{\theta}_0)\pi_0(\boldsymbol{\theta}_0)d\boldsymbol{\theta}_0} ,$$

where, under model $M_i$, $\boldsymbol{X}$ has density $f_i(\boldsymbol{x} \mid \boldsymbol{\theta}_i)$ and the unknown parameters $\boldsymbol{\theta}_i$ in $M_i$ are assigned a prior density $\pi_i(\boldsymbol{\theta}_i), i = 0, 1$. For given prior model probabilities $Pr(M_0)$ and $Pr(M_1) = 1 - Pr(M_0)$, the posterior probability of, say, $M_0$ is

$$(2.2) \qquad Pr(M_0 \mid \boldsymbol{x}) = \left[1 + B_{10} \frac{Pr(M_1)}{Pr(M_0)}\right]^{-1} .$$

In objective Bayesian analyses $\pi_i(\boldsymbol{\theta}_i)$ is chosen in an objective or conventional fashion and the hypotheses would be assumed to be equally likely a priori.

Use of objective priors has a long history in Bayesian inference (see, for example, [8] and [17] for justifications and references). They are, however, typically improper and are only defined up to an arbitrary multiplicative constant. This is not a problem in the posterior distribution, since the same constant appears in both the numerator and the denominator of Bayes theorem and so cancels. In model selection and hypothesis testing, however, it can be seen from (2.1) that when at least one of the priors $\pi_i(\boldsymbol{\theta}_i)$ is improper, the arbitrary constant does not cancel, so that the Bayes factor is then arbitrary and undefined. An important exception to this arises in invariant situations for parameters occurring in all of the models; Berger et al. [7] show that use of the (improper) right Haar invariant prior is then permissible.

One of the ways to address this difficulty is to try to directly "fix" the Bayes factor by appropriately choosing the multiplicative constant, as in [13]. Popular methods (the *intrinsic Bayes factor* [5] and the *fractional Bayes factor* [20]) for fixing this constant arise as a consequence of "training" the improper priors into proper priors based on part of the data or of the likelihood. We refer to Berger and Pericchi [6] for a review, references and comparisons. Another possibility is to directly derive appropriate "objective" but proper distributions $\pi_i(\boldsymbol{\theta}_i)$ to use in model selection; see [2] and [15] for methods and references. This is the approach taken in this paper (with a slight exception in Section 5).

## *2.2. Specification and justification of the objective priors*

Returning to the testing of the Poisson ($M_0$) *vs.* the ZIP ($M_1$) models, i.e., testing

$$(2.3) \qquad M_0: \boldsymbol{X} \sim f_0(\boldsymbol{x} \mid \lambda) \quad vs. \quad M_1: \boldsymbol{X} \sim f_1(\boldsymbol{x} \mid \lambda, p),$$

the key issue is the choice of the priors $\pi_0(\lambda)$ and $\pi_1(\lambda, p) = \pi_1(\lambda)\,\pi_1(p \mid \lambda)$.

A frequent simplifying procedure (both for subjective and objective methods) is to take $\pi_0(\lambda)$ equal to $\pi_1(\lambda)$, that is, to give the same prior to the parameters occurring in all models under consideration. This, however, may be inappropriate, since $\lambda$ might have entirely different meanings under model $M_0$ and under model $M_1$; the fact that we have used the same label does not imply that they have the same meanings. This frequent mistake is discussed, for example, in [7].

It has been argued that, if the common parameters are *orthogonal* to the remaining parameters in each model (that is, the Fisher information matrix is block diagonal), then they can be assigned the same prior distribution ([15], [16]). In this case, improper priors can be used, since the arbitrary constant would cancel in the Bayes factor.

Unfortunately, $p$ and $\lambda$ in the ZIP model are not orthogonal. We first reparameterize the original model. With $p^* = p + (1-p)e^{-\lambda}$, we rewrite $f_1(x \mid \lambda, p)$ as

$$(2.4) \qquad f_1^*(x \mid \lambda, p^*) = p^* I(x = 0) + (1 - p^*) f^T(x \mid \lambda), \quad x = 0, 1, 2, \ldots,$$

where $f^T(x \mid \lambda)$ is the zero-truncated Poisson distribution with parameter $\lambda$. Note that $p^* \geq e^{-\lambda}$. We can trivially express the Poisson ($M_0$) model as:

$$(2.5) \qquad f_0^*(x \mid \lambda) = e^{-\lambda} I(x = 0) + (1 - e^{-\lambda}) f^T(x \mid \lambda), \quad x = 0, 1, 2, \ldots,$$

and now it can intuitively be seen that $\lambda$ has the same meaning in both $f_1^*$ and $f_0^*$. Indeed the Fisher Information matrix for $p^*$ and $\lambda$ can be checked to be diagonal.

With an orthogonal reparameterization, Jeffreys (1961) recommended using (i) *Jeffreys prior* (the square root of Fisher information) for the "common" parameters; and (ii) a reasonable *proper* prior for the extra parameters in the more complex model.

The situation here is very unusual, however, in that the Jeffreys prior for the "common" $\lambda$ is different for each model. The *Jeffreys prior* for $\lambda$ in the Poisson model is well known to be $\pi_J^0 = 1/\sqrt{\lambda}$, whereas the Jeffreys prior for the orthogonalized ZIP model is easily shown to be the same as the Jeffreys prior for the truncated distribution $f^T(x \mid \lambda)$, which is

$$\pi_J^1(\lambda) = \frac{k(\lambda)}{\sqrt{\lambda}}, \quad \text{where} \quad k(\lambda) = \frac{\{1 - (\lambda + 1)e^{-\lambda}\}^{1/2}}{1 - e^{-\lambda}}.$$

That these priors are different after orthogonalization is highly unusual and can be traced to the fact that $\lambda$ also enters into the definition of the nested model, through $p^* = e^{-\lambda}$. In any case, we are left without clear guidance as to whether $\pi_J^0$ or $\pi_J^1$ should be used as the prior for $\lambda$. (Note that, in computing the Bayes factor, the same prior for $\lambda$ must be used in both the numerator and the denominator; otherwise one is facing the indeterminacy issues discussed earlier.)

Under the orthogonalized ZIP model, we also need to specify a proper prior for $p^*$ given $\lambda$, which we propose to take uniform over the interval $(e^{-\lambda}, 1)$, that is,

$$\pi_1(p^* \mid \lambda) = \frac{I(e^{-\lambda} < p^* \leq 1)}{1 - e^{-\lambda}}.$$

We can thus write the overall priors being considered for the two models $f_0^*(x \mid \lambda)$ and $f_1^*(x \mid \lambda, p^*)$ as, respectively,

$$\pi_0^l(\lambda) = \frac{k(\lambda)^l}{\sqrt{\lambda}}, \quad \pi_1^l(\lambda, p^*) = \frac{k(\lambda)^l}{\sqrt{\lambda}} \frac{I(e^{-\lambda} < p^* \leq 1)}{1 - e^{-\lambda}},$$

where $l$ is 0 or 1 as we utilize one or the other of the two Jeffreys priors for $\lambda$.

It is computationally more convenient to work in the original $(p, \lambda)$ parameterization. A change of variables above then results in the priors

(2.6) $$\pi_0^l(\lambda) = \frac{k(\lambda)^l}{\sqrt{\lambda}}, \quad \pi_1^l(\lambda, p) = \frac{k(\lambda)^l}{\sqrt{\lambda}} I(0 < p \leq 1),$$

which we will henceforth consider (for $l$ equal to 0 or 1).

We are not aware of any desiderata that would suggest a preference for either the $l = 0$ prior or the $l = 1$ prior, but luckily the two yield almost the same answers. Indeed, simple algebra shows that $k(\lambda)$ is a strictly increasing function of $\lambda$ and that

(2.7) $$\inf\ k(\lambda) = \frac{1}{\sqrt{2}} = 0.71 \quad \text{and} \quad \sup\ k(\lambda) = 1.$$

Thus $k(\lambda)$ is quite flat as a function of $\lambda$, so that $k(\lambda)^1$ and $k(\lambda)^0 = 1$ are very similar. An immediate consequence for the Bayes factors $B_{10}^l$, $l = 0, 1$ is that

$$B_{10}^0/\sqrt{2} \leq B_{10}^1 \leq \sqrt{2}\, B_{10}^0,$$

so that the two Bayes factors can only differ by a modest amount (and in practice the difference is much smaller than this).

It is obviously a bit simpler to work with the $l = 0$ prior, so we drop the $l$ superscript and henceforth utilize the prior

(2.8) $$\pi_0(\lambda) = \frac{1}{\sqrt{\lambda}}, \quad \pi_1(p, \lambda) = \frac{1}{\sqrt{\lambda}} I(0 < p \leq 1).$$

### 2.3. Objective Bayes factor for Poisson versus ZIP models

Recall that the model $M_0$ is the standard Poisson model and the model $M_1$ is the ZIP model. For a sample of $n$ counts $X_1, \ldots, X_n$, let $\boldsymbol{X}$ denote the sample, $k = \sum_{i=1}^n I(X_i = 0)$ be the number of zero counts, and $s = \sum_{i=1}^n X_i$ be the total count. Note that $k = n$ is equivalent to $s = 0$. For given data $\boldsymbol{x}$, the densities $f_0(\boldsymbol{x} \mid \lambda)$ and $f_1(\boldsymbol{x} \mid \lambda, p)$ under the two models are given by

$$f_0(\boldsymbol{x} \mid \lambda) = \frac{e^{-n\lambda}\lambda^s}{\prod_{i=1}^n x_i!}, \quad f_1(\boldsymbol{x} \mid \lambda, p) = \frac{[p + (1-p)e^{-\lambda}]^k (1-p)^{n-k} e^{-(n-k)\lambda} \lambda^s}{\prod_{i=1}^n x_i!}.$$

For $s > 0$ (i.e., the counts are not all zero),

$$m_0(\boldsymbol{x}) = \int f_0(\boldsymbol{x} \mid \lambda) \pi_0(\lambda) d\lambda = \frac{\Gamma(s + \frac{1}{2})}{n^{s+\frac{1}{2}} \prod x_i!}.$$

Using the binomial expansion of $[p + (1-p)e^{-\lambda}]^k$,

$$
\begin{aligned}
m_1(\boldsymbol{x}) &= \int f_1(\boldsymbol{x} \mid \lambda, p)\pi_1(p, \lambda) dp\, d\lambda \\
&= \frac{1}{\prod x_i!} \sum_{j=0}^{k} \frac{k!}{j!(k-j)!} \int_0^\infty \int_0^1 p^j(1-p)^{n-j} e^{-(n-j)\lambda} \lambda^{s-\frac{1}{2}} dp\, d\lambda \\
&= \frac{k!}{(n+1)! \prod x_i!} \sum_{j=0}^{k} \frac{(n-j)!}{(k-j)!} \Gamma(s+\tfrac{1}{2})(n-j)^{-(s+\frac{1}{2})}.
\end{aligned}
$$

Both $m_0(\boldsymbol{x})$ and $m_1(\boldsymbol{x})$ are finite and the Bayes factor $B_{10}(\boldsymbol{x}) = m_1(\boldsymbol{x})/m_0(\boldsymbol{x})$ is

$$
(2.9) \qquad B_{10}(\boldsymbol{x}) = \frac{k!}{(n+1)!} \sum_{j=0}^{k} \frac{(n-j)!}{(k-j)!} \left(1 - \frac{j}{n}\right)^{-(s+1/2)}.
$$

Note that, as intuitively expected, for any given $n$ the Bayes factor is increasing in $s$ (total count) for any fixed $k$ (the number of zero's), and is increasing in $k$ for any fixed $s$. We use (2.9) to calculate the Bayes factors for the examples in Section 3.

When $s = 0$ or equivalently all counts are zero ($\boldsymbol{x} = \boldsymbol{0}$), there is a problem. While $m_0(\boldsymbol{0}) = \Gamma(1/2)/\sqrt{n}$ remains finite, it is easy to see that $m_1(\boldsymbol{0})$ is infinite. Indeed for *any* prior of the form $h(p)\pi(\lambda)$, where $\pi(\lambda)$ is improper and $h(p)$ is a proper density (as is required for testing), the marginal density $m_1(\boldsymbol{0})$ will be infinite. This is because, for $\boldsymbol{x} = \boldsymbol{0}$, the density $f_1(\boldsymbol{x} \mid \lambda, p) \geq p^n$ implying $m_1(\boldsymbol{0}) \geq \int_0^1 p^n h(p) dp \int_0^\infty \pi(\lambda) d\lambda = \infty$. We discuss what to do for this case in Section 5.

## 3. Applications

In this section we apply our methodology to two datasets to detect if zero-inflation is present in the data. These examples have been analyzed for zero-inflation previously using both frequentist and Bayesian procedures. Since there are non-zero counts in both examples, the Bayes factors are computed using (2.9).

**Example 3.1.** The first dataset is the Urinary Tract Infection (UTI) data used in Broek [10], which used a score test to detect zero-inflation in a Poisson model. The data are collected from 98 HIV-infected men treated at the Department of Internal Medicine at the Utrecht University hospital. The number of times they had a urinary tract infection was recorded as $X$. The data are recorded in Table 1. Merely by looking at the data it is apparent that zero-inflation is present.

Equation (2.9) yields a Bayes factor $B_{10} = 223.13$ in favor of model $M_1$ versus model $M_0$; if the models were believed to be equally likely a priori, the resulting posterior model probabilities would be $Pr(M_1 \mid \boldsymbol{x}) = 0.995$ and $Pr(M_0 \mid \boldsymbol{x}) = 0.005$. This is indeed strong evidence in favor of the ZIP model.

In Bayesian testing of $\mathcal{H}_0 : p \leq 0$ versus $\mathcal{H}_1 : p > 0$, Bhattacharya et al. [9] obtained $Pr(p > 0 \mid \boldsymbol{x}) = 0.999$. The observed value of the score statistic was reported as 15.34 [10], yielding a $p$-value of 0.0001. All three analyses present strong evidence in favor of the ZIP model, but notice that the *p*-value seems to suggest stronger evidence against the Poisson null than the Bayesian analysis, and the point null Bayesian analysis suggests weaker evidence than the interval Bayesian test.

TABLE 1
*UTI Data*

| $X$ | 0 | 1 | 2 | 3 | Total |
|---|---|---|---|---|---|
| Frequency | 81 | 9 | 7 | 1 | 98 |

TABLE 2
*Terror Data*

| $X$ | 0 | 1 | 2 | 3 | 4 | Total |
|---|---|---|---|---|---|---|
| Frequency | 38 | 26 | 8 | 2 | 1 | 75 |

**Example 3.2.** The next dataset we consider is the Terrorism data from [11]. Table 2 gives the number of incidents of international terrorism per month ($X$) in the United States between 1968 and 1974. It is not intuitively clear whether or not there is zero-inflation in this data set.

The Bayes factor here is $B_{10} = 0.28$, yielding an objective posterior probability $Pr(M_1 \mid \boldsymbol{x}) = 0.219$, which actually supports the Poisson model. A previous analysis found $Pr(p > 0 \mid \boldsymbol{x}) = 0.507$, an indeterminate value [9]. The observed value of the score statistic is 0.04, with a *p*-value of 0.83. Conigliani et al. [11] test a Poisson null model against a nonparametric alternative, finding a fractional Bayes factor $B_{10}^F$ of 0.0089 of the nonparametric alternative to the Poisson; the apparent strength of this conclusion, compared with the other results, is rather puzzling.

## 4. Model selection in ZIP regression

Many applications involve count data where covariate information is available; see, for example, [14] and [18]. In this section we consider selecting between Poisson regression and ZIP regression models given by

$$M_0^R: \quad X_i \stackrel{ind}{\sim} Poisson(\lambda_i), \ i = 1, \ldots, n, \tag{4.1}$$

$$M_1^R: \quad X_i \stackrel{ind}{\sim} ZIP(\lambda_i, p), \ i = 1, \ldots, n. \tag{4.2}$$

For a known offset variable $a_{0i}$, a $q \times 1$ vector of covariates $\boldsymbol{a}_i$ and regression parameters $\boldsymbol{\beta} = (\beta_1, \ldots, \beta_q)^T$, suppose the $\lambda_i$ follow the log-linear relationship

$$\log(\lambda_i) = a_{0i} + \boldsymbol{a}_i^T \boldsymbol{\beta}.$$

We assume that the matrix $\boldsymbol{A}^T = (\boldsymbol{a}_1, \ldots, \boldsymbol{a}_n)$ is of rank $q$. Let $k$ denote the number of zero counts in the data. For simplicity of notation, we index the observations in such a way that all the zeros are given by the first $k$ counts.

### 4.1. Objective priors for model selection

Generalizing the argument in Section 2.2 to the regression case is easy in one case, but difficult in the other. If we choose to base the analysis on the Jeffreys prior for $\boldsymbol{\beta}$ under the Poisson regression model $M_0^R$, the generalization is straightforward: the Jeffreys prior is easily computed as

$$\pi_0^R(\boldsymbol{\beta}) = |\sum_{i=1}^n \lambda_i \boldsymbol{a}_i \boldsymbol{a}_i^T|^{1/2}. \tag{4.3}$$

Note that this prior is positive since the rank of $\boldsymbol{A}$ is $q$. Also, utilizing this prior for $\boldsymbol{\beta}$ under model $M_1^R$, along with the independent uniform prior for $p$, results in the following priors to be utilized to compute $B_{10}$:

$$(4.4) \quad \pi_0^0(\boldsymbol{\beta}) = |\sum_{i=1}^n \lambda_i \boldsymbol{a}_i \boldsymbol{a}_i^T|^{1/2}, \quad \pi_1^0(\boldsymbol{\beta},p) = |\sum_{i=1}^n \lambda_i \boldsymbol{a}_i \boldsymbol{a}_i^T|^{1/2} I(0 < p \leq 1).$$

The generalization to the regression case of the second prior considered in Section 2.2 is much more difficult, because the Jeffreys prior under the ZIP regression model is very complicated. In Section 2.2, the derivation of the corresponding Jeffreys prior was essentially done by ignoring the zero counts, utilizing only the truncated Poisson distribution. This suggests modifying (4.3) by removing the terms corresponding to the zero counts, resulting in

$$(4.5) \quad \pi_1^R(\boldsymbol{\beta}) = |\sum_{i=k+1}^n \lambda_i \boldsymbol{a}_i \boldsymbol{a}_i^T|^{1/2}.$$

From another intuitive perspective, the zero counts arising from the inflation factor are clearly irrelevant in fitting the log linear model to the $\lambda_i$ and, since we do not know which zero counts arise from the inflation factor, dropping them all from the Jeffreys prior has an appeal. Let $\boldsymbol{A}_+ = (\boldsymbol{a}_{k+1}, \ldots, \boldsymbol{a}_n)^T$. The prior (4.5) can only be used provided it is positive, which is ensured if the rank of $\boldsymbol{A}_+$ is $q$.

The resulting overall prior for use in computing $B_{10}$ is then

$$(4.6) \quad \pi_0^1(\boldsymbol{\beta}) = |\sum_{i=k+1}^n \lambda_i \boldsymbol{a}_i \boldsymbol{a}_i^T|^{1/2}, \quad \pi_1^1(\boldsymbol{\beta},p) = |\sum_{i=k+1}^n \lambda_i \boldsymbol{a}_i \boldsymbol{a}_i^T|^{1/2} I(0 < p \leq 1).$$

The first basic issue in use of these priors is whether or not they yield finite marginal distributions. This is addressed in the following theorems, the first of which deals with the marginal density under the Poisson regression model.

**Theorem 4.1.** *For the Poisson regression model and either the Jeffreys prior ($j = 0$) or the modified Jeffreys prior ($j = 1$),*

$$(4.7) \quad m_0^R(\boldsymbol{x}) = \int_{R^q} \prod_{i=1}^n \{\frac{e^{-\lambda_i}\lambda_i^{x_i}}{x_i!}\} \pi_j^R(\boldsymbol{\beta}) d\boldsymbol{\beta} < \infty.$$

*Proof.* See the Appendix. □

Note that with more than one covariate there is typically no closed-form expression for $m_0^R(\boldsymbol{x})$. Hence $m_0^R(\boldsymbol{x})$ needs to be evaluated by numerical or Monte Carlo integration.

For the ZIP regression model, the marginal density $m_1^R(\boldsymbol{x})$, under an arbitrary improper prior $\pi(\boldsymbol{\beta})$ for $\boldsymbol{\beta}$ and an independent uniform prior for $p$, is given by

$$(4.8) \quad m_1^R(\boldsymbol{x}) = \int_{R^q} \int_0^1 f_1(\boldsymbol{x} \mid \boldsymbol{\beta}, p) \pi(\boldsymbol{\beta}) \, dp \, d\boldsymbol{\beta},$$

where the density of $\boldsymbol{x}$, under model $M_1^R$, is given by

$$f_1(\boldsymbol{x} \mid \boldsymbol{\beta}, p) = \prod_{i=1}^k \{p + (1-p)e^{-\lambda_i}\}(1-p)^{n-k} \prod_{i=k+1}^n \frac{e^{-\lambda_i}\lambda_i^{x_i}}{x_i!}.$$

Again, as for $m_0^R(\boldsymbol{x})$, there is usually no closed-form expression for $m_1^R(\boldsymbol{x})$ and the marginal needs to be computed via numerical or Monte Carlo integration.

To investigate the finiteness of $m_1^R(\boldsymbol{x})$, note first that

$$(4.9) \qquad p^k(1-p)^{n-k} \prod_{i=k+1}^{n} \frac{e^{-\lambda_i} \lambda_i^{x_i}}{x_i!} \leq f_1(\boldsymbol{x} \mid \boldsymbol{\beta}, p) \leq \prod_{i=k+1}^{n} \frac{e^{-\lambda_i} \lambda_i^{x_i}}{x_i!} .$$

In view of this inequality and the independent uniform prior for $p$, the marginal $m_1^R(\boldsymbol{x})$ is finite if and only if

$$(4.10) \qquad \int_{R^q} \prod_{i=k+1}^{n} \frac{e^{-\lambda_i} \lambda_i^{x_i}}{x_i!} \pi(\boldsymbol{\beta}) \, d\boldsymbol{\beta} < \infty .$$

Theorem 4.2 below gives sufficient conditions for this to be finite under the priors (4.3) and (4.5) respectively. Recall that the $k$ zeros in the sample are labeled to correspond to the first $k$ observations. A key condition will be that the matrix $\boldsymbol{A}_+$ has rank $q$ which implies that $n \geq k + q$ (analogous to the condition of at least one positive count for the case of no covariate treated in Section 2).

**Theorem 4.2.** *Using $\pi_0^R(\boldsymbol{\beta})$: Suppose that, for the observation $X_j, j = 1, \ldots, k$, corresponding to the zero counts, the corresponding covariate vector $\boldsymbol{a}_j$ is such that*

$$(4.11) \qquad \boldsymbol{a}_j = \sum_{m=k+1}^{n} c_{mj} \, \boldsymbol{a}_m \quad \text{with} \quad c_{mj} \geq 0, j = 1, \ldots, k, m = k+1, \ldots, n.$$

*Then the marginal $m_1^R(\boldsymbol{x})$ is finite.*
*Using $\pi_1^R(\boldsymbol{\beta})$: If $\boldsymbol{A}_+$ has rank $q$, the marginal $m_1^R(\boldsymbol{x})$ is finite.*

*Proof.* See the Appendix. □

Clearly the condition under which $m_1^R(\boldsymbol{x})$ is finite is more general and much easier to check for $\pi_1^R(\boldsymbol{\beta})$ than for $\pi_0^R(\boldsymbol{\beta})$. This, together with the intuitive appeal of $\pi_1^R(\boldsymbol{\beta})$, leads us to recommend its use in practice. (Note that either of the two priors reduces to the prior recommended in Section 2 for the non-regression case.)

**Remark 4.1.** If the condition (4.11) fails, the marginal density $m_1^R(\boldsymbol{x})$ based on the Jeffreys prior may be infinite. For example, consider $n = 3$ and $q = 2$, with $\lambda_1 = \lambda_2^{c_1} \lambda_3^{c_2}$, $\lambda_2 = \exp(\beta_1)$, $\lambda_3 = \exp(\beta_2)$ for suitable nonzero $c_1, c_2$ to be chosen later. Then the determinant of information matrix for $\boldsymbol{\beta}$ is given by

$$|\boldsymbol{I}(\boldsymbol{\beta})| = \lambda_2 \lambda_3 + c_1^2 \lambda_2^{c_1} \lambda_3^{c_2+1} + c_2^2 \lambda_2^{c_1+1} \lambda_3^{c_2} ,$$

so that $|\boldsymbol{I}(\boldsymbol{\beta})|^{1/2} \geq |c_1| \lambda_2^{c_1/2} \lambda_3^{(c_2+1)/2}$. If $X_1 = 0$, $X_2 = x_2$ and $X_3 = x_3$, then

$$\begin{aligned} m_1^R(\boldsymbol{x}) &\geq \frac{|c_1|}{2} \int_{R^2} \frac{e^{-\lambda_2} \lambda_2^{x_2}}{x_2!} \frac{e^{-\lambda_3} \lambda_3^{x_3}}{x_3!} \lambda_2^{c_1/2} \lambda_3^{(c_2+1)/2} d\boldsymbol{\beta} \\ &= \frac{|c_1|}{x_2! x_3! 2} \int_0^\infty e^{-\lambda_2} \lambda_2^{x_2-1+.5c_1} d\lambda_2 \int_0^\infty e^{-\lambda_3} \lambda_3^{x_3-1+.5c_2+.5} d\lambda_3 = \infty , \end{aligned}$$

providing that $x_2 \leq -.5c_1$ or that $x_3 \leq -.5 - .5c_2$. For example, if $c_1 = -5$ and a sample produces $x_2 = 2$, then $m_1^R(\boldsymbol{x}) = \infty$. Note that here $\boldsymbol{a}_1 = -5\boldsymbol{a}_2 + c_2\boldsymbol{a}_3$, with $\boldsymbol{a}_2 = (1, 0)^T$ and $\boldsymbol{a}_3 = (0, 1)^T$, so that the condition (4.11) does not hold.

*4.2. An illustrative application*

We apply the methodology recommended in Section 4.1 to a dataset involving the number of AIDS-related deaths in men. The data provides the number of deaths for 598 census tracts in a large city of Spain over a period of eight years. The dataset, which was supplied to us by Dr. M.A.M. Beneyto, has a large number of tracts with zero deaths (actually, 303, which is $k$ in our notation). Along with the number of deaths, the dataset also provides, for each census tract, the expected number of deaths $E$ from AIDS (adjusting for the population and the distribution of ages in each tract) and an auxiliary variable $W$ (continuous in nature) measuring the social status of each census tract.

In our application and for the $i$th census tract, we take $\log(E_i)$ as the offset $a_{0i}$ and propose a log-linear regression for $\lambda_i$ with $q = 2$ and $\boldsymbol{a}_i = (1, W_i)^T$. First, we will ignore the covariate $W$ and compute the Bayes factor taking $q = 1$ and $\boldsymbol{a}_i = 1$ based on the Jeffreys prior. This model modifies the common mean model of Section 2.2 by incorporating the offset variable in the mean, which is here given by $E_i \lambda$ with $\lambda = \beta_1$. The marginal $m_1(\boldsymbol{x})$ is computed by one-dimensional numerical integration. Although it has a closed-form expression, it is rather complicated and omitted here to save space. This expression is given in the Appendix in [1]. For the specific data here, $B_{10} = 22,975$ which gives overwhelming evidence in favor of the ZIP model.

Epidemiologists who are knowledgeable about this study believed that the large number of zero counts in the data could be explained by the covariate measuring the social status and, indeed, suspected that a ZIP regression model would not be needed if the covariate were incorporated into the analysis. The Bayes factor in favor of the ZIP regression model versus the Poisson regression model (with $q = 2$) is given by 7.25. While this Bayes factor provides a moderate amount of evidence in favor of the ZIP regression model, it is much smaller than 22,975, indicating that, indeed, the covariate can explain most of the excess zero counts.

In this example, it is possible that the same inflation parameter $p$ may not be appropriate for all individuals. Just like using the log-linear models for $\lambda_i$, we can treat each $p_i$ differently (as $p$ may change according to the covariates) and fit a logistic regression model for $p_i$. But it is highly likely that there would be severe confounding between the two regressions, which is particularly problematical with objective Bayesian analysis (since there is not a proper subjective prior to overcome the confounding).

## 5. Analysis with insufficient positive counts

As noted in Section 2, the marginal density under model $M_1$ based on an improper prior for $\lambda$ is not finite when all counts are zeros, and hence the Bayes factor is not well-defined. This is not a difficulty of only model selection; in this situation, it is also not possible to make inferences about the parameters of the ZIP model, since the joint posterior of the parameters (under the ZIP model) is improper. Indeed, when all counts are zero, the ZIP model parameters are not identifiable, and the data do not provide enough information to estimate the parameters. Since objective Bayes methods are typically based on information from the data alone, it is not surprising that problems are encountered.

We could simply invoke this argument and refrain from considering the case when all counts are zero. However, it is interesting to explore several methodologies that have been proposed for difficult testing situations, partly to judge the success of the methodologies and partly to try to provide a reasonable answer to this case. We continue, throughout the section, to assume that $p \sim Un(0,1)$.

## 5.1. All zero counts in the non-regression case

We mentioned that to resolve the identifiability issue in the ZIP model for the data with all zeros we need a proper prior on $\lambda$. This can be done by either subjectively specifying a proper prior for $\lambda$ or by "training" the improper priors into proper priors based on part of the data or of the likelihood. In particular, the intrinsic Bayes factor approach [5] utilizes a part of the data as a training sample to train the improper prior to get a proper posterior. Although this approach works successfully in many examples, it is not successful in the present problem. Our investigation of this approach [1] is omitted here to save space. We discuss below the case where a subjective proper prior on $\lambda$ is specified based on certain considerations.

If a proper prior is needed to define the Bayes factor for the situation of all zero counts, the most direct approach is to find a proper prior that seems compatible with certain behaviors that we expect of the Bayes factor in this situation. A natural proper prior to consider for $\lambda$ is a Gamma ($Ga(a,b)$) conjugate prior under the Poisson model ($M_0$) given by the Gamma $g(\lambda \mid a, b)$ density

$$g(\lambda \mid a, b) = \frac{b^a e^{-b\lambda} \lambda^{a-1}}{\Gamma(a)},$$

where $a, b$ are suitably chosen positive constants. Of course, one is welcome to simply make subjective choices here, but we will argue for a certain choice (or choices) based on rather neutral thinking.

First, we assume that the *same* gamma prior is appropriate for $\lambda$, both under the Poisson and the ZIP models. This can be justified by the orthogonalization argument used in Section 2.2. With the uniform density for $p$ and the $Ga(a,b)$ prior for $\lambda$, the resulting Bayes factor for arbitrary data $\boldsymbol{x}$ can be computed to be

$$(5.1) \quad B_{10}(\boldsymbol{x}) = \frac{k!}{(n+1)!} \sum_{j=0}^{k} \frac{(n-j)!}{(k-j)!} \left(1 - \frac{j}{n+b}\right)^{-(s+a)},$$

by a similar argument to that leading to (2.9). This Bayes factor includes as a special case the objective Bayes factor in (2.9); indeed the Jeffreys prior used there was a limiting case of the $g(\lambda \mid a, b)$ for $a = 1/2$ and $b = 0$. Note that the Bayes factor (5.1) is increasing in $s$, $k$ and $a$, and decreasing in $b$.

For the special case $\boldsymbol{x} = \boldsymbol{0}$ (that is $s = 0$ and $k = n$), note that $f_1(\boldsymbol{0}|\lambda, p) \geq f_0(\boldsymbol{0}|\lambda)$. Hence, using the same proper prior for $\lambda$ with both the Poisson and the ZIP models, it follows that $m_1(\boldsymbol{0}) \geq m_0(\boldsymbol{0})$, and hence, $B_{10}(\boldsymbol{0}) \geq 1$. In particular, for the $Un(0,1)$ prior for $p$ and $Ga(a,b)$ prior for $\lambda$, it can be checked that

$$(5.2) \quad B_{10}(\boldsymbol{0}) = \frac{(n+b)^a}{n+1} \sum_{j=0}^{n} \frac{1}{(j+b)^a} \geq 1.$$

This is reasonable: when a long stream of *only* zeros is observed, it is entirely natural to say that the data favor the ZIP model. But the degree of favoritism depends on $a$ and $b$, and we turn to rather speculative desiderata to narrow the choice. Recall that the mean of the $Ga(a,b)$ distribution for $\lambda$ is $ab^{-1}$ and the variance is $ab^{-2}$.

In order for the prior not to be too sharp, it is reasonable to require the prior standard deviation to be no less than the prior mean. This implies that $a \leq 1$. It also seems reasonable to require the prior mean to be at least 1, so that small values of $\lambda$ do not have excessive prior probability. This leads to $b \leq a$. Since the Bayes factor is decreasing in $b$, the smallest Bayes factor satisfying the above constraints (that is, the one lending the most support for the Poisson model $M_0$) is then obtained by taking $b = a$ (this gives a prior mean of 1). It is not unreasonable to select this prior as it belongs to a reasonable class which is most favorable to the null model. Finally, one might judge it to be unappealing to utilize a prior for $\lambda$ which is not bounded near zero (for $a < 1$ the gamma density is decreasing with an asymptote at $\lambda = 0$) which implies that $a$ should be at least 1. Thus we end up with the choice $a = b = 1$. Note that $a = 1$ is the upper limit of $a \leq 1$ and the choice $a = 1$ now counterbalances the Bayes factor in favor of $M_1$ (whereas $b = a$ in the range $b \leq a$ tilts the Bayes factor in favor of $M_0$). This reasoning is all rather speculative and, of course, the result is a particular prior, which may not reflect actual prior beliefs. Nevertheless it is instructive to study the behavior of the Bayes factor when this prior is used.

For $a = b = 1$, that is, the Exponential(1) distribution, it can be checked that $B_{10} = \sum_{j=0}^{n}(j+1)^{-1}$, which is thus our recommended default Bayes factor when observing only zero counts. Note that $B_{10}(\mathbf{0}) \approx \log(n+1)$ for large $n$. So a *large* string of all zero counts in a sample will lead to a Bayes factor approaching infinity at the slow rate of $\log(n)$. The large sample behavior of the Bayes factor for this type of sample seems intuitively reasonable.

### *5.2. Insufficient positive counts in the regression case*

In the regression situation of Section 4, it was necessary to have sufficient positive counts so that the conditions of Theorem 4.2 were satisfied. We will restrict discussion here to the situation involving the prior specifications in (4.6), for which the key condition needed for the marginal to be finite was that the matrix $\mathbf{A}_{+}((n-k) \times q)$ should be of rank $q$. If the number of positive counts $n - k$ is insufficient so that $t$, the rank of $\mathbf{A}_{+}$, is less than $q$, this solution will not work.

**Remark 5.1.** Indeed, neither the prior for $\boldsymbol{\beta}$ given by (4.3) nor by (4.5) guarantees a finite positive marginal density. We omit the proof to save space. A proof may be found in the Appendix in [1].

We call this situation one of rank deficiency, with the rank deficiency of $\mathbf{A}_{+}$ equal to $q - t$. The situation is analogous to the case of all zero counts without covariates discussed in Subsection 5.1. (In the setup of that section, $q = 1$ and rank $\mathbf{A}_{+}$ less than 1 means that $k = n$, i.e., no positive counts.) We could again merely recognize that this type of data is just not informative enough to allow for objective Bayes analysis. We shall however propose a prior that yields finite marginal densities, following similar reasoning to that used in Section 5.1.

We continue to use a $Un(0,1)$ prior for $p$ and focus on proposing suitable priors for $\boldsymbol{\beta}$. A discussion similar to that in subsection 5.1 shows that this prior has to be at least, partially proper.

Note that, instead of specifying a prior on $\boldsymbol{\beta}$, we can specify a prior on $q$ independent parametric functions of $\boldsymbol{\beta}$; our specific proposal is to carefully choose these functions such that $t$ of them are well identified by the data with positive counts while the remaining $q - t$ are not. We then propose to use a version of Jeffreys prior on the former $t$ functions, and a proper prior on the latter $q - t$ functions.

Specifically, let $\boldsymbol{A}_0$ denote the $k \times q$ matrix whose $k$ rows are $\boldsymbol{a}_1^T, \ldots, \boldsymbol{a}_k^T$. A rank of $\boldsymbol{A} = q$ and a rank of $\boldsymbol{A}_+ = t$ imply a rank of $\boldsymbol{A}_0 \geq q - t$. Let $V_+ \subseteq R^q$ denote the vector space of dimension $t$ formed by the columns of $\boldsymbol{A}_+^T$. Suppose $\boldsymbol{a}_{i_1}, \ldots, \boldsymbol{a}_{i_r}$ are all of the vectors from $\boldsymbol{a}_1, \ldots, \boldsymbol{a}_k$ corresponding to the zero counts which are in $V_+$. Note that $0 \leq r \leq k - (q - t)$. These vectors are linear combinations of the vectors $\boldsymbol{a}_{j_1}, \ldots, \boldsymbol{a}_{j_t}$ and the corresponding $\lambda_{i_1}, \ldots, \lambda_{i_r}$ are functions of $\lambda_{j_1}, \ldots, \lambda_{j_t}$. From the set of $\{\lambda_j : j \in \{1, \ldots, k\} - \{i_1, \ldots, i_r\}\}$ we select $q - t$ $\lambda$'s, $\lambda_{l_1}, \ldots, \lambda_{l_{q-t}}$ such that $\{\boldsymbol{a}_{j_1}, \ldots, \boldsymbol{a}_{j_t}, \boldsymbol{a}_{l_1}, \ldots, \boldsymbol{a}_{l_{q-t}}\}$ is linearly independent.

Note that there is an $(n - k) \times t$ matrix $\boldsymbol{C}$ of rank $t$ such that

$$(\boldsymbol{a}_{k+1}, \ldots, \boldsymbol{a}_n) = (\boldsymbol{a}_{j_1}, \ldots, \boldsymbol{a}_{j_t})\boldsymbol{C}^T.$$

Let $\boldsymbol{D} \equiv \boldsymbol{D}(\lambda_{j_1}, \ldots, \lambda_{j_t})$. Then, the information matrix for $\lambda_{j_1}, \ldots, \lambda_{j_t}$ based on the Poisson model for the observations $k+1, \ldots, n$ is given by

$$(5.3) \qquad \boldsymbol{I}(\lambda_{j_1}, \ldots, \lambda_{j_t}) = \boldsymbol{D}^{-1} \boldsymbol{C}^T Diag(\lambda_{k+1}, \ldots, \lambda_n) \boldsymbol{C} \boldsymbol{D}^{-1}.$$

We define a partial Jeffreys prior for $\lambda_{j_1}, \ldots, \lambda_{j_t}$ by

$$(5.4) \qquad \pi_{PJ}(\lambda_{j_1}, \ldots, \lambda_{j_t}) = \{\prod_{i=1}^{t} \lambda_{j_i}^{-1}\} |\boldsymbol{C}^T Diag(\lambda_{k+1}, \ldots, \lambda_n) \boldsymbol{C}|^{1/2}.$$

Let $\{\boldsymbol{b}_1, \ldots, \boldsymbol{b}_{q-t}\}$ denote an orthonormal basis of the space spanned by $\boldsymbol{a}_{l_1}, \ldots, \boldsymbol{a}_{l_{q-t}}$. Define $\xi_w = e^{\boldsymbol{b}_w^T \boldsymbol{\beta}}$, $w = 1, \ldots, q - t$. Note that $\lambda_{l_w}, w = 1, \ldots, q - t$ can be expressed in terms of $\xi_1, \ldots, \xi_{q-t}$. Indeed,

$$\log(\lambda_{l_w}) = a_{0l_w} + \sum_{h=1}^{q-t} d_{wh} \log(\xi_h), \quad w = 1, \ldots, q - t,$$

where $d_{wh} = \boldsymbol{b}_h^T \boldsymbol{a}_{l_w}$. Finally, we assign independent exponential distributions with mean 1 to each of $\xi_1, \ldots, \xi_{q-t}$. This prior will induce a proper distribution on $\lambda_{l_w}, w = 1, \ldots, q - t$ with a density which we denote by $\pi_{prop}(\lambda_{l_1}, \ldots, \lambda_{l_{q-t}})$. The final prior used to calculate the marginal density under model $M_1^R$ is then given by

$$\pi(\lambda_{j_1}, \ldots, \lambda_{j_t}, \lambda_{l_1}, \ldots, \lambda_{l_{q-t}}) = \pi_{PJ}(\lambda_{j_1}, \ldots, \lambda_{j_t}) \pi_{prop}(\lambda_{l_1}, \ldots, \lambda_{l_{q-t}});$$

this is partially Jeffreys prior and partially proper. The corresponding prior density on $\boldsymbol{\beta}$ is, of course, obtained through transformation. Further, along the line of the proof of Theorem 4.2, it can be checked that the marginal density $m_1^R(\boldsymbol{x})$ will be finite. We omit the details to save space.

While there is arbitrariness in the specific choice of $\lambda_{l_1}, \ldots, \lambda_{l_{q-t}}$ to assign a subjective prior distribution based on exponential distributions, the partial Jeffreys prior in (5.4) remains invariant to the choice of $t$ independent $\lambda$'s from $\lambda_{k+1}, \ldots, \lambda_n$. This solution thus seems reasonable for small $q - t$.

To avoid the arbitrariness, we could consider all possible selections of $(q - t)$ of the $\lambda$'s from $\lambda_1, \ldots, \lambda_k$ so that these $q - t$ and $t$ of the $\lambda$'s from $\lambda_{k+1}, \ldots, \lambda_n$ define a reparameterization of $\boldsymbol{\beta}$. For each selection we can calculate the Bayes factor, and in the spirit of IBF we can take a suitable average over all these Bayes factors. If the rank deficiency of $\boldsymbol{A}_+$ is 1, we will have $k - r$ Bayes factors to average.

**Appendix**

*Proof of Theorem 4.1.* From (4.3) and (4.5) it is immediate that $\pi_1^R(\boldsymbol{\beta}) \leq \pi_0^R(\boldsymbol{\beta})$. Thus it is enough to prove (4.7) for $j = 0$. Let $\boldsymbol{i}$ denote the indices $(i_1, \ldots, i_q)$ and $\boldsymbol{A}(\boldsymbol{i})$ denote a $q \times q$ submatrix of $\boldsymbol{A}$ based on rows $i_1, \ldots, i_q$. Then by Binet-Cauchy expansion of determinant (cf. Noble [19], p. 226) it can be shown that

$$\text{(A1)} \qquad |\sum_{i=1}^n \lambda_i \boldsymbol{a}_i \boldsymbol{a}_i^T| = \sum (\lambda_{i_1} \ldots \lambda_{i_q})|\boldsymbol{A}(\boldsymbol{i})\boldsymbol{A}(\boldsymbol{i})^T|,$$

where the summation is over all submatrices of order $q \times q$. Dropping the terms from the above summation for which $|\boldsymbol{A}(\boldsymbol{i})\boldsymbol{A}(\boldsymbol{i})^T| = 0$ we get from (4.3) that

$$\text{(A2)} \qquad \pi_0^R(\boldsymbol{\beta}) \leq \sum{}^{*} (\lambda_{i_1} \ldots \lambda_{i_q})^{1/2} |\boldsymbol{A}(\boldsymbol{i})\boldsymbol{A}(\boldsymbol{i})^T|^{1/2},$$

where $\sum^*$ denotes summation over all $q \times q$ matrices for which $|\boldsymbol{A}(\boldsymbol{i})\boldsymbol{A}(\boldsymbol{i})^T| > 0$.
Since $e^{-\lambda_i} \lambda_i^{x_i}/x_i! < 1$, from (4.7) and (A2) we get

$$\text{(A3)} \qquad m_0^R(\boldsymbol{x}) \leq \sum{}^{*} \int_{R^q} \prod_{j=1}^q \{\frac{e^{-\lambda_{i_j}} \lambda_{i_j}^{x_{i_j}}}{x_{i_j}!}\} (\lambda_{i_1} \ldots \lambda_{i_q})^{1/2} |\boldsymbol{A}(\boldsymbol{i})\boldsymbol{A}(\boldsymbol{i})^T|^{1/2} d\boldsymbol{\beta}.$$

Recall that $\log(\lambda_i) = a_{0i} + \boldsymbol{a}_i^T \boldsymbol{\beta}$. Now transforming $\boldsymbol{\beta}$ to $(\lambda_{i_1}, \ldots, \lambda_{i_q})$ and using the Jacobian of transformation $(\lambda_{i_1} \ldots \lambda_{i_q})^{-1} |\boldsymbol{A}(\boldsymbol{i})\boldsymbol{A}(\boldsymbol{i})^T|^{-1/2}$, we get from (A3) that

$$\text{(A4)} \qquad m_0^R(\boldsymbol{x}) \leq \sum{}^{*} \prod_{j=1}^q \int_0^\infty \frac{e^{-\lambda_{i_j}} \lambda_{i_j}^{x_{i_j}-.5}}{x_{i_j}!} d\lambda_{i_j} < \infty,$$

since each of the integrals in the right hand side of (A4) is finite. This completes the proof of Theorem 4.1. □

*Proof of Theorem 4.2.* First, as in (A1) and (A2), it can be shown that for some positive $c$ (not depending on parameters) less than 1

$$\text{(A5)} \qquad \begin{aligned} c \sum{}^{*} (\lambda_{i_1} \ldots \lambda_{i_q})^{1/2} |\boldsymbol{A}(\boldsymbol{i})\boldsymbol{A}(\boldsymbol{i})^T|^{1/2} \\ \leq \pi_0^R(\boldsymbol{\beta}) \leq \sum{}^{*} (\lambda_{i_1} \ldots \lambda_{i_q})^{1/2} |\boldsymbol{A}(\boldsymbol{i})\boldsymbol{A}(\boldsymbol{i})^T|^{1/2}. \end{aligned}$$

In view of this inequality and (4.10), the marginal $m_1^R(\boldsymbol{x})$ is finite if and only if

$$\text{(A6)} \qquad \int_{R^q} \prod_{i=k+1}^n \frac{e^{-\lambda_i} \lambda_i^{x_i}}{x_i!} (\lambda_{i_1} \ldots \lambda_{i_q})^{1/2} |\boldsymbol{A}(\boldsymbol{i})\boldsymbol{A}(\boldsymbol{i})^T|^{1/2} d\boldsymbol{\beta} < \infty$$

for each $\boldsymbol{i} = (i_1, \ldots, i_q)$ for which $|\boldsymbol{A}(\boldsymbol{i})\boldsymbol{A}(\boldsymbol{i})^T| > 0$.

Note that the sufficient condition stated in the theorem and the condition that rank of $\boldsymbol{A}$ is $q$ imply that the regression matrix $\boldsymbol{A}_+^T = (\boldsymbol{a}_{k+1}, \ldots, \boldsymbol{a}_n)$ corresponding to the set of positive counts has rank $q$.

Suppose, with no loss of generality, $i_1 < \cdots < i_q$ in (A6). Also, suppose $i_1 < \cdots < i_u \leq k < i_{u+1} < \cdots < i_q$. It is possible that $u$ may be 0 or may be $q$.

By the assumed condition that for $j = 1, \ldots, k$, $\boldsymbol{a}_j$ can be expressed as a linear combination of $\boldsymbol{a}_{k+1}, \ldots, \boldsymbol{a}_n$ with nonnegative coefficients, it follows that

$$\lambda_{i_j} = h_{i_j} \prod_{m=k+1}^{n} \lambda_m^{c_{m i_j}}, \quad j = 1, \ldots, u,$$

where $c_{m i_j} \geq 0$ and $h_{i_j} > 0$. Then

$$\prod_{j=1}^{u} \lambda_{i_j} = f \prod_{m=k+1}^{n} \lambda_m^{b_m},$$

where $b_m = \sum_{j=1}^{u} c_{m i_j} \geq 0$ and $f > 0$ are free from parameters.

Then the integrand (without $|\boldsymbol{A}(\boldsymbol{i})\boldsymbol{A}(\boldsymbol{i})^T|^{1/2}$) in (A6) can be simplified as

$$\prod_{i=k+1}^{n} \frac{e^{-\lambda_i} \lambda_i^{x_i}}{x_i!} (\lambda_{i_1} \ldots \lambda_{i_q})^{1/2}$$

$$= \prod_{i=k+1}^{n} \frac{e^{-\lambda_i} \lambda_i^{x_i + \frac{1}{2} b_i}}{x_i!} (\lambda_{i_{u+1}} \ldots \lambda_{i_q})^{1/2}$$

$$(A7) \qquad = \left[ \prod_{j=u+1}^{q} \frac{e^{-\lambda_{i_j}} \lambda_{i_j}^{x_{i_j} + \frac{1}{2} b_{i_j} + \frac{1}{2}}}{x_{i_j}!} \right] \left[ \prod_{l=1}^{n+u-k-q} \frac{e^{-\lambda_{\alpha_l}} \lambda_{\alpha_l}^{x_{\alpha_l} + \frac{1}{2} b_{\alpha_l}}}{x_{\alpha_l}!} \right],$$

where $\{\alpha_1, \ldots, \alpha_{n+u-k-q}\} = \{k+1, \ldots, n\} - \{i_{u+1}, \ldots, i_q\}$.

Suppose $\{s_1, \ldots, s_q\} \subset \{k+1, \ldots, n\}$ is such that $\{\boldsymbol{a}_{s_1}, \ldots, \boldsymbol{a}_{s_q}\}$ is a linearly independent set (such a set exists since $\boldsymbol{A}_+$ is of rank $q$). Note that for $y > 0$ the function $g(u) = e^{-u} u^y$ is maximized at $u = y$ implying

$$(A8) \qquad e^{-u} u^y \leq e^{-y} y^y \quad \text{for all } u > 0.$$

By (A8) we get from (A7) that

$$(A9) \qquad \prod_{i=k+1}^{n} \frac{e^{-\lambda_i} \lambda_i^{x_i}}{x_i!} (\lambda_{i_1} \ldots \lambda_{i_q})^{1/2} \leq D \left( \prod_{j=1}^{q} e^{-\lambda_{s_j}} \lambda_{s_j}^{d_{s_j}} \right),$$

where $D > 0$ is a constant independent of the parameters and $d_{s_j} = x_{s_j} + \frac{1}{2} b_{s_j} + \frac{1}{2}$ if $s_j \in \{i_{u+1}, \ldots, i_q\}$, and $d_{s_j} = x_{s_j} + \frac{1}{2} b_{s_j}$ if $s_j \in \{\alpha_1, \ldots, \alpha_{n+u-k-q}\}$.

The Jacobian of transformation from $\boldsymbol{\beta}$ to $\lambda_{s_1}, \ldots, \lambda_{s_q}$ is $H/(\lambda_{s_1} \ldots \lambda_{s_q})$ for some $H > 0$ constant. Then since $d_{s_j} \geq 1$ for $j = 1, \ldots, q$, by (A9) we have

$$(A10) \qquad \int_{R^q} \prod_{i=k+1}^{n} \frac{e^{-\lambda_i} \lambda_i^{x_i}}{x_i!} (\lambda_{i_1} \ldots \lambda_{i_q})^{1/2} d\boldsymbol{\beta} \leq HD \prod_{j=1}^{q} \int_0^{\infty} e^{-\lambda_{s_j}} \lambda_{s_j}^{d_{s_j}-1} d\lambda_{s_j} < \infty.$$

By (A10) and (A6) we conclude that $m_1^R(\boldsymbol{x})$ corresponding to $\pi_0^R(\boldsymbol{\beta})$ is finite. To prove finiteness of $m_1^R(\boldsymbol{x})$ corresponding to $\pi_1^R(\boldsymbol{\beta})$ note that by (4.10)

$$m_1^R(\boldsymbol{x}) \leq \int_{R^q} \left( \prod_{i=k+1}^{n} \frac{e^{-\lambda_i} \lambda_i^{x_i}}{x_i!} \right) \pi_1^R(\boldsymbol{\beta}) d\boldsymbol{\beta}.$$

Finiteness of the right hand quantity in the last display follows from a version of Theorem 4.1 corresponding to the prior $\pi_0^R(\boldsymbol{\beta})$ by replacing $n$ observations from the Poisson by $n - k$ observations from Poisson. This completes the proof. □

**Acknowledgments.** The authors would like to thank the Community of Valencia Group in the Project *Desigualdades socioeconómicas y medioambientales en ciudades en España, Proyecto MEDEA* for the data used in Section 4, a referee for valuable comments, and Archan Bhattacharya for computing help. Part of this research was conducted while Bayarri and Datta were visiting SAMSI/Duke University, whose support is gratefully acknowledged.